# Further examples of apodictic discourse, I

Satyanad Kichenassamy

**Abstract.** The analysis of problematic mathematical texts, particularly from India, has required the introduction of a new category of rigorous discourse, *apodictic discourse*. We briefly recall why this introduction was necessary. We then show that this form of discourse is widespread among scholars, even in contemporary Mathematics, in India and elsewhere. It is in India a natural outgrowth of the emphasis on non-written communication, combined with the need for freedom of thought. New results in this first part include the following: (i) Āryabhaṭa proposed a geometric derivation of a basic algebraic identity; (ii) Brahmagupta proposed an original argument for the irrationality of quadratic surds on the basis of his results on the *varga-prakṛti* problem, thereby justifying his change in the definition of the word *karaṇī*.

## 1. Issues and background.

This paper is a first of a series based on a course[1] given in Paris during the first semester of academic year 2021-22. It showed that the study of Indian Mathematics opened the way to a truly global history of Mathematics, by providing conceptual tools applicable to texts in all languages. It also spelled out the objective reasons why these tools have not been identified earlier: we build on a considerable body of research in Indology and Historical Epistemology; an examination of the historiography of Indian Mathematics will show the circumstances that explain this slow development. The lectures were in part based on papers that have appeared recently, but also contain some new material and new perspectives on earlier results. We focus here on this material, giving some background information along the way to make the paper self-contained. The focus of interest is the existence of several different strategies and forms of discourse adapted to the communication of rigorous mathematical results and their justification. Two extreme forms may be distinguished. One is *dogmatic discourse,* rigorous, conclusive and unmotivated, the lack of motivation being made up by imposing upon the reader a particular sequence of arguments. At the other end of the spectrum is *apodictic discourse,* that is also rigorous, conclusive but motivated. Because it presents results in a natural order, it suffices to make the logical organization of the discourse clear, pointing out difficult steps and essential conditions, because a well-prepared reader will be able to fill in the details without ambiguity. Apodictic discourse may therefore be significantly shorter than dogmatic discourse, making use of discursive strategies to encode much necessary information in the statements of the results, rather than through an auxiliary discourse[2]. Since mathematical texts may be rigorous without being dogmatic[3] in this sense, it is necessary to introduce a second form of discourse beside dogmatic discourse[4].

---

[1] This course on "Apodictic discourse and the Historiography of Mathematics" was delivered during the first semester of academic year 2021-22 in Paris, within Jan. E.M. Houben's chair on "Sources and history of the Sanskrit tradition" at the École Pratique des Hautes Études (University Paris Sciences et Lettres). Most of the following is based on lectures delivered at the Sorbonne in January and February 2022, but also includes material from talks at the Harish-Chandra Institute, as part of their 75[th] Independence celebration, and at the "International Conference on Ancient Indian Astronomy and Mathematics" organized by SVTB College, Palakkad, as well as lectures delivered in December 2019 at the Chennai Mathematical Institute. Thanks to Jan Houben, C.S. Seshadri, Parvathy K.P. and Chandan Dalawat for their kind invitations.
[2] It therefore contains its *apodeixis* (the conclusive argument of a proof in Greek) within itself, hence its name.
[3] The term should be taken in its etymological sense (< *dokein*); it need not have negative connotations. Both forms of discourse have their uses.
[4] The first example of this type of discourse analysis is "Baudhāyana's rule for the quadrature of the circle", *Historia Mathematica*, **33**(2) (2006) 149–183. The term 'apodictic discourse' was introduced in *Aestimatio*, **13** (2016–2018), 119–140, reprinted as "Translating Sanskrit mathematical texts" with a biographical notice in *Aestimatio N.S.,* **1** (2020), 183-204.

Our work builds on what happened in the last fifty years or so in History of Science or research-level Mathematics, as well as in Indology. We now understand that much of what used to be uncritically taken for granted in Science and Epistemology ultimately rests on discursive constructions of collective or personal identities [5]. Since constructions of identity are smokescreens that hamper direct contact with the texts, leading to anachronistic readings[6], we focus in this first part on several new examples showing how rigorous Indian texts actually work, spelling out the concepts that are operative in them. Only later shall we show how the rediscovery of ancient texts of all cultures in the nineteenth century made older representations of Mathematics untenable, and made it clear that there were in India several "epistemic cultures" with different concepts and different demands for rigor.

We are only interested here in the communication, through discourses, of mathematical results at the highest level of rigor, it being understood that heuristic reasoning and, indeed, unconscious reasoning that is never verbalized, play a major role in mathematical invention[7]. They are normally suppressed in rigorous discourse, especially when the audience and the author have a common mathematical and cultural background. Rigor may also be relaxed in situations where audiences are merely interested in directions for further action, not with justifications: a patient expects prescriptions from his or her physician, not a lecture on physiology or anatomy. However, the physician needs a rigorous discourse because he or she knows that anything less may lead to serious errors. In an Indian context, the same base texts may be used by various audiences, which may create the illusion that the crudest understanding is necessarily correct, since some practitioners care for no other. However, the texts we consider are so carefully written, and their results so complex that one cannot assume that the derivations have been left implicit. If they are not given out as separate arguments, we must conclude that they are encoded in the discursive structure – in the order and choice of words, or through the deliberate insertion of anomalies that should make the reader think. While these discursive strategies are very common in Indian Philosophy, Grammar or Literary analysis, we show that they also should be identified in other fields of scholarship, particularly Mathematics. We are therefore interested in the professional mathematician that knows that anything less than inference (*anumāna*) – be it perceptual evidence (*pratyakṣa*), arguments of analogy (*upamāna*) or arguments based on the authority of tradition (*āgama*) – may simply be wrong. Any discourse derives its validity from its coherence, as expressed in a pattern of inference. Rigorous communication relies on the analysis of other people's discourses, and on the construction of one's own discourses.

---

[5] On the discursive construction of identities, see for instance Stuart Hall, *Race, Ethnicity, Nation,* Cambridge U. Press, 2017 and, of course, E.W. Said, *Orientalism*, Routledge, London, 1978. In modern science, there are several "epistemic cultures" that form subgroups within a wider scientific community (Karin Knorr-Cetina, *Epistemic Cultures: How the sciences make knowledge,* Harvard University Press, Cambridge, MA, 1999; Evelyn Fox Keller, *Making Sense of Life: Explaining Biological Development with Models, Metaphors, and Machines,* Harvard University Press, Cambridge, MA, 2002). This notion seems to provide a finer description than Kuhn's "normal science". A similar phenomenon was pointed out by Sarasvati Amma, about Indian Geometry and (*Geometry in Ancient and Medieval India*, second revised ed.: Motilal Banarsidass, New Delhi, 1999, p. 81) and, later, by the author (in *Historia Mathematica,* in 2010 and 2012 for India, and in 2015 for Italy). These schools may fail to understand one another. Similarly, in Hellenistic Mathematics, the problems solved by Diophantos do not seem to be related to earlier mathematical texts in Greek.

[6] This problem was identified only recently by Niccolò Guicciardini (*Anachronism in the History of Mathematics,* Cambridge U. Press, 2021, see his introduction). In the present paper, anachronism will always be taken in its negative sense: the attribution to an author of knowledge or notions that he or she could not possibly have possessed or used. We suggest that anachronism is closely related to the discursive construction of identities, and that, for the examples we consider, anachronism is part of a *fantasy-echo* in the sense of Joan W. Scott.

[7] The basic study here is Jacques Hadamard's monograph (*The Psychology of Invention in the Mathematical Field*, Dover, New York, 1954). It reflects not only his experience, but also that of well-known mathematicians that took part in his survey. Similar observations have been made by Gustave Choquet ("Les processus mentaux de la creation", *Séminaire de Philosophie et Mathématiques,* fasc. 4, (1994) 1–21.). The conclusions are easy to check by any mathematician, on the basis of his or her own professional experience.

The prototype of an apodictic discourse is Brahmagupta's[8] elaborate discourse on the cyclic quadrilateral, in his *Brāhmasphuṭasiddhānta* [9] (from 628), Prop. 12.21-32. Despite what has been claimed for several centuries, even in India, his formulation is not faulty: he did specify that his quadrilaterals had to be cyclic, by introducing a neologism, *tricaturbhuja* "triquadrilateral[10]", that refers to a trilateral completed into a quadrilateral by adding a fourth point on its circumcircle, as is clear from Prop. 12.27[11]. The structure of this proposition is such that there is no ambiguity as to the meaning of the term. But it was misunderstood by a hostile commentator in the ninth century and its meaning forgotten and with it, the elaborate derivation he suggested. In view of the length of his discourse, we focus in this part on shorter examples that also serve to introduce basic tools in Brahmagupta's and, more generally, Indian Mathematics. The fact that Brahmagupta's mathematics was misunderstood and indeed, unduly criticized by later authors is proof enough that a loosely defined common cultural background is not enough to create an epistemic culture and therefore, that essentialism does not provide an adequate description of the facts.

Cultural background may nonetheless bear on scientific developments[12]. While apodictic discourse is found in many cultures, its emergence in India seems to be intimately related to familiar features of Indian scholarship[13], such as the asymmetric relation between student and master (*guru*) and the emphasis on verbal, unwritten communication. Indeed, to express geometric theorems without writing or drawing, it is necessary to replace concrete objects by abstract ones, and to carry out reasoning directly on mental objects[14]. And to maintain freedom of thought, arguments should be inferred by the student, by the free exercise of his or her intelligence in scrutinizing the master's words, who remains available to dispel doubts if needed. In other words, the student will enjoy freedom if it is granted by the *guru*.

---

[8] His full name is *Jiṣṇusutabrahmagupta* "Brahmagupta, son of Jiṣṇu". Like most Indian authors, Brahmagupta does not claim credit for his inventions, only for his exposition, as is indicated by his use of the first person. How much is due to his father is difficult to assess, although the mention of his father's name in the text suggests that he was a well-known scholar, and that Brahmagupta learnt Mathematics and Astronomy from him.

[9] Henceforth abbreviated as *BSS*. For the text and some background information, see H. T. Colebrooke, *Algebra, with Arithmetic and Mensuration, from the Sanskrit of Brahmegupta and Bhascara*, London: J. Murray, 1817; Sudhākara Dvivedin, *Brāhmasphuṭa Siddhānta and Dhyānagrahopadeśādhyāya of Brahmagupta*, Benares, 1902; R. S. Sharma (Chief Ed.), *Brāhma-Sphuṭa Siddhānta, with Vāsanā, Vijñāna and Hindi Commentaries*, New Delhi: Indian Institute of Astronomical and Sanskrit Research, 1966.

[10] We need a neologism in English to translate Brahmagupta's Sanskrit neologism, since no exact modern equivalent exists.

[11] Cf. S. Kichenassamy, "Brahmagupta's Derivation of the Area of a Cyclic Quadrilateral", *Historia Mathematica*, **37** (1) (2010) 28-61 and "Brahmagupta's propositions on the perpendiculars of cyclic quadrilaterals", *Historia Mathematica*, **39** (4) (2012) 387-404. Brahmagupta's many results and derivations go well beyond his area formula. They will be recalled in a later part of this paper, with complements.

[12] A striking recent example is Aldo G. Gargani's analysis of Newton's physics (*La science sans fondements : La conduite intellectuelle comme structuration de l'expérience commune*, (tr. C. Alunni) Vrin, Paris, 2013). For older examples of the interference of non-scientific issues in scientific work, see for instance J.J. Wunenburger (ed.) *Bachelard et l'épistémologie française*, Presses Universitaires de France, Paris, 2003 ; D. Lecourt, *L'épistémologie historique de Gaston Bachelard*, Vrin, Paris, 2002 or L. Daston, *Classical Probability in the Enlightenment*, Princeton University Press, 1995 among many others.

[13] The relation between science and culture has been recently revisited since the realization that much work in this direction aimed at establishing an ideological foundation for imperialism, by defining exclusive identities to replace much more fluid concepts of identity that existed all over the world, in India in particular. For an attempt to put the problem on a new basis, see the recent volume edited by K. Chemla and E. Fox-Keller (*Cultures without culturalism : The making of scientific knowledge,* Duke University Press, Durham and London, 2017). This volume clearly outlines the issues even though Indian science is not considered in it. To this line of thought, we add that earlier work neglected the existence of a *reflexive* dimension in ancient India that provides new tools to tackle the problem, because the concept of the knowing subject, that may leave society and reflect on it from the outside, was problematized there very early. This dimension is missed if one focuses on "practical" aspects of science, interesting as they are, because of their dependence on special features of social structure. For related reasons, the notion of practice itself is being problematized (Theodore R.Schatzki, K. Knorr-Cetina and Eike von Schatzki (eds.), *The practice turn in contemporary theory,* Routledge, London and New York, 2001).

We give in this series of papers new examples and revisit older examples of apodictic discourse. The outline of this work is as follows. Section 2 recalls some terminology and gives two new examples of apodictic discourse, one from Ramanujan's works and the other, from Harish-Chandra's. Section 3 shows that the term *vyavahāra*, when it refers to sections of a mathematical work, refers to (apodictic) discourses, and that this meaning is consistent with other sources, while other translations are inconsistent. We then analyze, in Section 4, Āryabhaṭa's presentation of algebraic identities and show that it contains clear indications of his intended derivations. Section 5 shows how Brahmagupta, in the eighteenth chapter of BSS, relates the *varga-prakṛti* equation[15] $Nx^2 + k = y^2$ to the irrationality of certain square roots. Recall[16] that this eighteenth chapter introduces a new technical meaning for the term *karaṇī*, that refers for Baudhāyana to any square root: for Brahmagupta, it is limited to irrational square roots, for which he stresses that they are constructible, as perpendiculars of isosceles trilaterals. Brahmagupta related their irrationality to the existence of infinitely many solutions of the *varga-prakṛti* equation when $k = 1$. A conclusion summarizes the new elements in this first part.

## 2. Rigorous mathematics as discourse.

### 2.1. Terminology.

We spell out the definitions of some basic terms in our discussion to prevent any possible confusion with related terms, or uses in other contexts.

An *apodictic discourse* is a rigorous, conclusive, and motivated discourse. It may be contrasted with a *dogmatic discourse*, namely a rigorous, conclusive discourse, but unmotivated, or buttressed by artificial motivation. A discourse may be dogmatic for good reason: if the audience does not have sufficient background to appreciate the motivation, if they are merely interested in the results rather than their derivation, or for other reasons – for example, if the original motivation was an incorrect induction. An example is given by Ramanujan's paper on "round numbers"[17]. He simply lists earlier work, without motivation, and goes on describing his contribution. It is Hardy who gives his original motivation: this paper stemmed from an attempt to substantiate an imprecise argument that both Ramanujan and he had made independently, and which turned out to be misleading[18]. It is natural that Ramanujan

---

[14] Note that this view is not an idealism since private mental objects produced by one mind are not "essences" that are directly accessible to other minds.

[15] Here, *N* and *k* are given, and the problem is to find *x* and *y*.

[16] On this point, see the analysis of Prop. 18.38 of BSS in §5.2 of S. Kichenassamy's paper "Brahmagupta's apodictic discourse," *Gaṇita Bhāratī*, **41**(1) (2019), 93-113.

[17] "Highly composite numbers", *Proceedings of the London Mathematical Society*, vol. 2, XIV, (1915), 347-409; *Collected Papers of Srinivasa Ramanujan* (G.H. Hardy, P.V. Seshu Aiyar and B.M. Wilson eds., Cambridge Univ. Press, 1927), Paper n° 15, pp. 78-128.

[18] See G.H. Hardy (*Ramanujan. Twelve Lectures on Subjects suggested by his Life and Work*, Cambridge Univ. Press, 1940), lecture III, pp. 48-57, on round numbers; a round number is defined as "the product of a considerably large number of comparatively small factors" (p. 48). The motivation was this: "It is a matter of common observation that round numbers are very rare; the fact may be verified by anyone who will make a habit of factorising numbers, such as numbers of motor cars of railway carriages, which are presented to his attention in a random manner. Both Ramanujan and I had observed this phenomenon, which seems at first a little paradoxical." (*ibid.*). Hardy adds in a footnote (quotation marks are in the original): "Half the numbers are divisible by 2, …one-

should have adopted the dogmatic form while focusing on his results, suppressing a conjecture that turned out to be false, while Hardy preferred the apodictic form to put the result in context and thereby illustrate Ramanujan's *ingenium*.

Not every text is a discourse. A composite text cannot be treated as a single discourse. A code, a list, a compilation are not discourses although their organization may be significant and imply a discourse. A grammatical rule or an organized list of such rules, a mathematical formula or table, or a computer program are not discourses, even if it is possible to construct a discourse that explains their meaning or use. By contrast, a treatise on grammar; an article containing formulae; the user's manual for a software program; or the elaborate sets of comments found within some computer programs, may be analyzed as discourses. An apodictic discourse is the expression of the coherent view of one author (or of a collective acting as one), aimed at a specific audience. The analysis of apodictic discourse naturally takes on an abstract character because it is rests on inference (*anumāna*) on the basis of word-representations, the latter being mental constructs and therefore not accessible to sense perception (*pratyakṣa*).

*Dialogues* constitute an intermediate form, commonly found in commentaries of philosophical texts, where the master typically helps a student in *removing doubts*. This interaction shows that the student is not expected to follow the master blindly. However, such discourses may or may not be conclusive, depending on whether the student, real or imagined, has put forward all logically possible objections. The same goes for dialogues with opponents, in which the eristic dimension may be significant. There are many intermediate forms that need not be recalled.

A *proposition* is a mathematical statement, that could be a theorem or a problem. It could also include an indication, or contain a reference to earlier propositions etc.

A *derivation* is a motivated sequence of arguments from assumptions to conclusion. It is *rigorous* or *conclusive* if the sequence is complete. Not every derivation is conclusive. Thus, arguments relying on sensory perception are not conclusive, even if they may suggest the idea of a proper derivation.

A *proof* is an unmotivated sequence of arguments from assumptions to conclusion. It is *rigorous* or *complete* if the sequence of steps is unbroken, but for small gaps that are easy to fill[19]. An incomplete proof or derivation could be *heuristic* if it suggests a way to a complete one or *flawed* if at least one of its steps is incorrect and beyond repair. The absence of a derivation in a discourse does not always imply that it is encoded in the discursive structure. The derivation could be common knowledge. In some contexts, 'proof' and 'derivation' may be interchangeable if no confusion ensues.

There are also forms of non-rigorous discourses, such as *suggestive discourses*, that indicate a general line of argument. A suggestive discourse may contain the kernel of a rigorous argument and may then be called *heuristic*. But it may also be essentially flawed and is then *fallacious*. Not every text is a discourse: a set of notes, a table or a figure are not discourses; however, their full meaning may be disclosed to others only through explanations namely, through a discourse.

---

sixth by both 2 and 3. Surely then we may then expect most numbers to have a *large* number of factors? But the facts seem to show the opposite."

[19] It is not uncommon to leave gaps in proofs, so long as the audience can easily reconstruct the argument. See D. Fallis, "Intentional gaps in mathematical proofs", *Synthese* **134** (2003), 45-69, L.E. Andersen, "Acceptable gaps in mathematical proofs", *Synthese* **197**(1) (2020), 233-247.

## 2.2. An example of discourse analysis in contemporary Mathematics.

"Mathematics is not just a body of results, each one attached to a technical argument, but an intricate system of theories and a historical process"[20]. However, this system is not always apparent, to the extent that essential elements of mathematical discourse appear to be missing, either deliberately[21], or merely because "familiarity breeds contempt": theories that have been used for years become so hackneyed that they are hardly worth mentioning[22]. But there are situations in which the optimal form of mathematical communication, even when there is no desire to suppress intermediate steps, is not necessarily provided by a linear discourse in which every argument is spelled out. Here is a typical example, that describes a situation familiar to many of us. In a paper in Harish-Chandra's memory[23], Rebecca Herb[24] recalls how she encountered one of his famous works in the course of her thesis:

> I was lucky because my thesis advisor, Garth Warner, knew those papers well enough that instead of saying, "Read Harish-Chandra's papers," or "Read *Discrete series*, I[25]," he said, "Read p. 302 of *Discrete series,* I."

Observe that the advisor does not explain why this passage should be read, nor what would be gained by reading it. The passage indicated is to be found near the end of the paper, and it is not immediately clear why this would be a good place to start. Actually, the advisor gave the only piece of information that was not obvious if one wanted to start working on the problems he suggested. Herb continues:

> The good thing about Harish-Chandra's papers was that if you knew what you were looking for and where to start, everything was written down. […]. You might have to refer back to three or four of his earlier papers[26] for results, or even definitions, but he told you exactly where to look. I started on p. 302 of *Discrete series,* I with Lemma 56 and worked my way backward and forward, picking up a lemma here and there from earlier sections of the paper and from earlier papers. By the time I really understood that page, I was ready to write the first part of my thesis.[27]

Observe how all falls into place provided the student analyzes the paper and takes the logical steps suggested by her analysis: "everything was written down". This is strikingly similar to

---

[20] J. Gray, *Linear Differential Equations and Group Theory from Riemann to Poincaré*, 2nd ed., Birkhäuser, Boston, MA, (2000, repr. 2008), p. xviii.
[21] See D. Fallis, "Intentional gaps in mathematical proofs" (*Synthese* **134** (2003), 45–69) and W. Byers, *How Mathematicians Think: Using ambiguity, Contradiction, and Paradox to Create Mathematics* (Princeton University Press, Princeton, NJ, 2007, repr. 2010).
[22] The author has often heard Peter D. Lax say that "linearity breeds contempt": we are so familiar with features of linear partial differential equations (PDEs) that we are much more interested in those methods that also apply to nonlinear PDEs. To this day, there is no general theory of nonlinear partial differential equations comparable to the theory developed in the last century for linear PDEs with constant coefficients, and none is in sight. To a large extent, linear partial differential equations with variable, especially non-smooth coefficients, also cannot be treated in a unified way.
[23] Harish-Chandra हरिश्चन्द्र is one of the most important mathematicians of the twentieth century, Indian or otherwise. See Roger E. Howe's biographical memoir "Harish-Chandra, 1923-1983", National Academy of Sciences, Washington, D.C., 2011. http://www.nasonline.org/publications/biographical-memoirs/memoir-pdfs/harish-chandra.pdf
[24] Rebecca A. Herb, "Harish-Chandra and his work", *Bull. AMS* , **25** (1), (July 1991), 1-17, see page 1.
[25] Harish-Chandra, "Discrete series for semisimple Lie groups, I: Construction of invariant eigendistributions", *Acta Mathematica,* **113** (1965), 241-318.
[26] This refers to references [2, a–m] in *Discrete series,* I.
[27] *Ibid.*

Brahmagupta's writing: if one takes his work as a set of rules, it appears that essential information is missing. However, analysis shows that indeed, "everything" was there.

This work of analysis, similar to the close reading advocated by literary critics like G. Lanson, and taken up by many authors such as I. Richards and the Yale school, is no different from the careful reading of philosophical and literary texts often found in India[28].

### 2.3. Written or unwritten discourse?

Indian mathematical discourse has the special characteristic of being typically imparted orally, even when manuscripts were abundant[29], so that the student is encouraged to think without material support. This is conducive to the development of mental reasoning. The emphasis on unwritten transmission in Indian scholarship is familiar. We are interested in the implications of the avoidance of writing for the development of mathematical thought in India – even when writing material was abundant. F. Staal's essay[30] and its reviews may give an idea of earlier attempts at drawing the consequences of shunning writing. We depart from them in three respects; First, we limit ourselves to rational discourse, as opposed to hermeneutics. We deal with textual analysis, not interpretation. Second, following P.-S. Filliozat, the central issue is not orality[31] *per se*, but the primacy of unwritten transmission, in rational scholarly communication, over writing[32]. Third, we are interested in innovative thinking rather than faultless transmission. Thus, oral learning has several dimensions: training of memory to ensure the preservation of texts; developing mental reasoning; perceiving the many layers of knowledge in other peoples' discourses, especially those of one's own *guru*[33].

Independent thinking requires that the *guru* create a safe space for investigation. This may be achieved by refraining from spelling out the teaching. Here are examples where the guru does not spell out the instruction, and even withholds it when the student is not ready to receive it. First, in *Chāndogya Upaniṣad*, VII.1[34], "Nārada approached Sanatkumāra and said, 'Teach me, Venerable Sir,' He said, 'Come to me with (tell me) what you know. Then I will teach you what is beyond that.' [35]" Here, the guru tailors the level of the teaching to the student's level of knowledge. A second passage from the same *upaniṣad* shows that the master will not correct mistakes if student is not ripe for improvement: in VIII.7 sqq., Indra and Virocana, who sought out Prajāpati, both receive the same partial knowledge and leave. Only Indra returns for higher instruction, having perceived the need for it[36]. In a third passage (IV.9), Satyakāma similarly teaches his student fully only after he has reached a certain level and has perceived the need for knowledge. Thus, quite generally, knowledge should not be doled out before the desire for it (*jijñāsā*) arises. A fourth passage shows that mental reasoning does not necessarily lead to an idealistic view of knowledge. Vedic ritual has three correlated dimensions: action, language and silence (IV.16), represented by different actors. Therefore, language is not autonomous but is correlated with thinking and action. Knowledge, know-how and the representation of knowledge through language are neither independent nor separate[37].

---

[28] S. Kichenassamy, « L'Analyse Littéraire au service de l'Histoire des Mathématiques : Critique interne de la Géométrie de Brahmagupta », *Comptes-Rendus des Séances de l'Académie des Inscriptions et Belles-Lettres* (Communication du 20 avril 2012), **CRAI 2012**, II (avril-juin) (2012) p. 781-796.
[29] Parvathy K.P. stressed, in her concept note describing the project topic for the Palakkad conference that: "the significance of oral learning has been grossly ignored [in Kerala]; the result of this omission [in the educational system] has been disastrous."
[30] Frits Staal, *The Fidelity of Oral Tradition and the Origins of Science*, North-Holland, 1986.
[31] The term may be misleading because of its psychoanalytical dimension, present to any modern writer's mind.
[32] "Even if oral transmission is always appreciated, even if a composition in *sūtra* style and in verse is an aid to memorization, the pandits never refused writing, never neglected the help they could derive from it." (Pierre-Sylvain Filliozat, "Ancient Sanskrit mathematics : an oral tradition and a written literature", in *History of Science, History of Text*, ed. Karine Chemla, Springer, Dordrecht, 2004, pp. 137-157, see p. 148).
[33] Recall that the *guru* is, literally, the "man of weight" – *gurvī* for a woman – a person whose words carry weight. The relation of the student to the *guru* (and his wife) is strongly positive. To a large extent, one's *guru* is, in this system, the highest authority. Being granted freedom of thought by the *guru* is tantamount to complete freedom. This is what explains that innovative thinking should be compatible with the admission of the *guru* – and the *guru* alone – as supreme authority.
[34] The text is from S. Radhakrishnan, *The Principal Upaniṣads* (George Allen & Unwin, London, 1953 ; repr. Humanities Press, Atlantic Highlands, NJ, 1992). We follow him (and Śaṃkara) for the meaning.
[35] *adhīhi, bhagavaḥ, iti hopasasāda sanatkumāraṃ nāradaḥ, taṃ hovāca: yad vettha tena mopasīda, tatas ta ūrdhvam vakṣyāmīti, sa hovāca*

This general outlook has consequences for ancient mathematics. On the one hand, the avoidance of writing, and the need for a private analysis, perforce performed mentally, make it necessary to replace sensory objects and actions by word-representations, leading to the general validity of results, independently of any material representation. Hence the abstract and general nature of propositions found in Indian mathematical texts. Since the text is a combination of word-representations, an analysis of its coherence is naturally performed discursively as well. This leads to a further level of abstraction in which word-representations – abstracted from objects – are themselves manipulated by an intellectual process, rather than subjected to a formal process that could be carried out in writing. The emphasis on unwritten transmission does not mean that writing should be completely shunned: mathematical operations may be written, and some of Brahmagupta's words refer to writing[38].

Thus, focusing on the specific features of unwritten discourse helps clarify that sense perceptions need to be analyzed in order to produce meaning. In particular, visual evidence, such as a diagram, may be suggestive but is not conclusive unless its meaning is clarified by a discourse. This brings us back to the analysis of discourse. We next show that a familiar term in Indian Mathematics, *vyavahāra*, also refers to discourse.

### 3. Brahmagupta's *vyavahāra* as discourse

Chapter 12 of *Brāhmasphuṭasiddhānta* opens with a description of the topics covered in it, namely twenty *parikarma* and eight *vyavahāra*, defining the mathematician in the process.

12.1. *parikarma viṃśatiṃ yaḥ saṅkalitādyāṃ pṛthagvijānāti |*
*aṣṭau ca vyavahārān chāyāntān bhavati gaṇakaḥ saḥ ||*

Who, the twenty elementary operations[39] (*parikarma*)
Beginning with addition, severally comprehends,
As well as the eight *vyavahāra*[40]
Of which shadow is the last, is a mathematician. (12.1)

---

[36] This teaching concerns the existence and properties of the three states: waking, dreaming, and dreamless sleep.
[37] See IV.16.2: *tayor anyatarām manasā saṃskaroti brahmā, vācā hotā'dhvaryur udgītā anyatarām…* "Of these, the *Brahmā* performs one with his mind; by speech the *Hotṛ*, the *Adhvaryu* and the *Udgātṛ* the other." Thus, ritual activity is performed by four actors that represent mind, speech and action. These three dimensions are parallel. Of the four actors, the *Brahmā* is supposed to keep silent unless an error is committed. He then points it out and prescribes appropriate measures.
[38] One example among many: *adho'dho sthāpyam* "should be set down one under the other" in BSS 18.3-6.
[39] The elementary operations, described in I.2-13 are (a) addition, subtraction, multiplication, division, squaring, extraction of square roots, cubing, extraction of cube roots; (b) reduction of five types of complex fractions to the standard form $p/q$; (c) rule of three (direct and inverse), rules of five, seven, nine and eleven quantities, barter.
[40] Namely: *miśraka-, śreḍhī-, kṣetra-, khāta-, citi-, krākacika-, rāśi-* and *chāyā-vyavahāra*.

And indeed, these topics are considered in order by Brahmagupta, who completes them by a section of a more fundamental character, that provides a transition to other chapters. We focus here on the import of the term *vyavahāra*.

Each of the *vyavahāra* is a self-contained discourse on a particular mathematical topic. The word *vyavahāra* may be used in the mathematical literature to qualify results or methods fit for common use, but this meaning is inappropriate when *vyavahāra* is contrasted with *parikarman*. Now, the first meaning of *vyavahāra* is "verbal usage"[41]: it is that which is established by verbal intercourse – among scholars in this case –, and refers to a mental object, not a practice in the naive sense. Thus, "[a]ll everyday verbal usage (*vyavahāra*) takes place by means of word meaning (*padārtha*) that is dealt with as literally true (*mukhyeneva*), but is in fact a reality (*bhāva*) created by conceptual differentiation (*vikalpa*)."[42] Patañjali's Yogasūtra I.9 makes it clear that such mental constructs (*vikalpa*) are objects of language, that may be *vastuśūnyaḥ* (devoid of concrete reality) and, as Vyāsa makes it clear, are neither real nor unreal[43]. Thus, *vyavahāra* does not refer to "practice" but to "verbal usage", a conceptual representation that may not be mapped to an empirical object. In a scientific context, a *vyavahāra* is therefore a discourse that represents a conceptual model of a real situation; its rigor derives from its discursive consistency, not from its "practical" efficiency. It follows that a *vyavahārika* mathematical result is not merely a rough approximation, but an approximate value of which the degree of inadequacy may be ascertained by known processes of representation. Thus, rigorous plane geometry is only approximate as a model for an approximately spherical Earth. We can account for all mathematical meanings of *vyavahāra* in this manner.

The term *vyavahāra* is also quite common in business and legal contexts, in which the same meaning seems to lie at the root of all others. Indeed, Menski's study[44] of this term indicates that: "[a] necessary re-reading of the relevant ancient texts therefore suggests, in [his] view, that what we learn about *vyavahāra* is not a prescriptive model, or a set of rules for how to do legal business. […] Rather, then, these texts should be read as a guidance […] in how to deal with doubts over righteousness and appropriate forms of behavior for [Indians] in all life situations." The terms "intercourse", "business" or "procedure" similarly have related, complex semantic fields in English, but none of them exactly corresponds to the term *vyavahāra*.

All this suggests that *vyavahāra*, in a mathematical context, is the Sanskrit equivalent of "(rational) discourse (based on common concepts)": a rational treatment of a topic based on common scholarly knowledge of the times[45]. Its derivative *vyavahārika* may be rendered by "conventional" – or more explicitly, "the result of a consensus obtained discursively" – rather than "practical", since its opposites, theory or creation[46], are neither mentioned nor relevant. It thus appears that *vyavahāra*, for Brahmagupta, represented a form of discourse that presents a mathematical theory in the form of a coherent exposition. Similarly, his material on *parikarman* may be scrutinized as discourse that provides the conceptual basis for the following *vyavahāra*. Thus one may translate *kṣetravyavahāra* by "(apodictic) discourse on closed figures".

---

[41] J.E.M. Houben, *The* Sambandha-samuddeśa *(Chapter on relation) and Bhartṛhari's Philosophy of Language*, Egbert Forsten, Groningen, 1995, see p. 333 (n. 522), that also relies on the authority of Subramania Iyer. P.-S. Filliozat also confirmed to me that this meaning of the term is correct.
[42] *Kārikā* 82 (Houben, *op. cit*, p. 310).
[43] Vyāsa's gloss on Patañjali's Yogasūtra I.9 reads *śabdajñānānupātī vastuśūnyo vikalpaḥ* (The creation of the mind is the [psychological transformation] that follows the knowledge of words and is devoid of reality). Vyāsa's gloss further explains that this *vikalpa* is the basis of *vyavahāra,* even though it is neither provable nor contradicts the provable (*sa na pramāṇopārohi na viparyayopārohi ca |vastuśūnyatve'pi śabdajñānamāhātmyanibandhano vyavahāro dṛśyate…*) (The text and gloss are from P.-S. Filliozat, *Yogabhāṣya de Vyāsa sur le Yogasūtra de Patañjali,* Éditions Āgamāt, Palaiseau, 2005, pp. 52-55 ; we merely translated his French into English).
[44] Werner Menski, "On *vyavahāra*" (*Indologica Taurinensia*, **33** (2007), 123-147), see p. 144. This view is compatible with the usual etymological explanation of *vyavahāra* < *vi-ava-hṛ* "remove (a thorn etc.)", see p. 138.
[45] The term "topic" has a similar meaning if we take the term in the sense implied by Aristotle's opening lines of his *Topics,* (A 1, 100a 18-20) that read: "The purpose of the present treatise is to discover a method by which we shall be able to reason from generally accepted opinions about any problem set before us and shall ourselves, when sustaining an argument, avoid saying anything self-contradictory." Ἡ μὲν πρόθεσις τῆς πραγματείας μέθοδον εὑρεῖν ἀφ᾽ ἧς δυνησόμεθα συλλογίζεσθαι περὶ παντὸς τοῦ προτεθέντος προβλήματος ἐξ ἐνδόξων, καὶ αὐτοὶ λόγον ὑπέχοντες μηθὲν ἐροῦμεν ὑπεναντίον. (Text and translation from E.S. Forster, *Aristotle. Topics.* Loeb Classical

While Brahmagupta's *BSS*, more than nine times longer than the *Āryabhaṭīya*, makes it relatively easy to recognize the articulation of discursive structure and mathematical content, this is also possible in the latter work, as we show next. We focus on a passage that simply does not make sense as a record of a rules, but does if we consider that it suggests a derivation of the results it presents.

### 4. An example of apodictic discourse in the Āryabhaṭīya.

Among the basic tools of Indian mathematics after Āryabhaṭa are identities enabling transformations of two given quantities into two other quantities of the same kind. For Brahmagupta[47], the following are basic: first, if the sum and difference of two quantities are known, then these quantities are known (this is the *saṅkramaṇa* or "transition" from sum and difference to the two quantities); second, if the difference and the difference of squares of two quantities are known, then we can recover them (*viṣamakarman* or "asymmetric operation", since one of the given quantities involves squaring, not the other). The second builds on the former. Indeed, if $a^2 - b^2$ and $a - b$ are known, so is $(a^2 - b^2)/(a - b) = a + b$, so that we may recover $a$ and $b$ by *saṅkramaṇa*. The expansion of a square in the form $(u + v)^2 - u^2 = (2u + v)v = 2uv + v^2$ follows from $a^2 - b^2 = (a + b)(a - b)$, by applying it with $a = u + v$ and $b = u$. Let us see how these results come up in Āryabhaṭa's Prop. II.23 and II.24[48].

II.23.   *samparkasya hi vargād viśodhayedeva vargasamparkam |*
   *yattasya bhavatyardhaṃ vidyād guṇakārasaṃvargam ||*

   From the square of the sum
   One should indeed simply subtract the sum of the squares.
   That which is [its] half
   Should be known to be the product of the factors.

In modern symbols, this expresses the relation
$$ab = \frac{1}{2}[(a+b)^2 - (a^2 + b^2)].$$
But is this modernizing reduction a fair representation of the text? Let us consider its wording more closely. Note the terms for sum (*samparka*: 'join, concatenation') and product (*saṃvarga* 'arranged in classes' as in an array), that suggest a representation of two quantities $a$ and $b$ by lines, and their product by an array or an oblong (see Fig. 1). The formulation suggests a

---

Library, Harvard Univ. Press, 1960). Unfortunately, "topic" also has other meanings so that it would be confusing to equate *vyavahāra* with it.

[46] In Greek, *prattein* (action) may be contrasted with *theorein* (contemplation), but also with *poiein* (creation).
[47] See BSS 18.37 and our commentary on it in *Gaṇita Bhāratī* **41** (2019), 93-113.
[48] For the text, see K.S. Shukla and K.V. Sarma, *Āryabhaṭīya of Āryabhaṭa,* Indian National Science Academy, New Delhi, 1976.

geometric derivation. Āryabhaṭa had already stressed in Prop. II.3 that *varga* refers to a (geometric) square, to its area, as well as to the product of two equals. The emphasis (*hi* indeed, *eva* simply) suggests that this derivation is considered simpler than an alternative that would presumably rest on the distributivity of multiplication, that Brahmagupta stresses (*BSS* 12.56). Once the appropriate diagram has been identified (Fig. 1), the derivation is immediate: from a square of side $a + b$, it suffices indeed to subtract the 'concatenation' of the two squares – that have a point in common on the figure – to obtain two oblongs, each of which has area $ab$, hence the result. Let us turn to the next proposition.

II.24   *dvikṛtiguṇāt saṃvargāddvyantaravargeṇa saṃyutānmūlam |*
   *antarayuktaṃ hīnaṃ tadguṇakāradvayaṃ dalitam ‖*

> Multiplying the square of two by the product,
> (And) adding the square of the difference of two (quantities), the root (of the result)
> Increased or decreased by (their) difference,
> [Yields] this pair of multipliers by halving.

In modern symbols, this expresses the relation:

$$a, b = \frac{1}{2}\left\{\sqrt{2^2\, ab + (a-b)^2} \pm (a-b)\right\}.$$

The question is: what relation to his earlier result is Āryabhaṭa suggesting here? There is a striking anomaly: four is expressed as the square of two, which sounds quaint indeed[49]. The metre by itself does not impose it: one could have replaced[50] *dvikṛtigu(ṇāt)* par *kṛtaguṇitāt* or *vedaguṇāt*, which replaces *dvikṛtigu* a *gaṇa* of *sarvalaghu* type, by another one, or by a *gaṇa* of the form *gll*, both admissible for an *ārya* metre[51]. Therefore, the indirect expression of 4 as square of 2 is a deliberate choice that contains an information to be discovered by analysis[52].

Now, if we remember that the *saṅkramaṇa* and *viṣamakarman* were basic tools of Indian mathematics at the time, it is easy to derive proposition II.24. For the result implies that the radical is equal to the sum $a + b$ and therefore, $2^2 ab + (a-b)^2 = (a+b)^2$. Now, from the *viṣamakarman*, the difference of squares $(a+b)^2 - (a-b)^2 = (a+b+a-b)(a+b-a-b) = (2a) \times (2b) = 2^2 ab$. Thus, the mention of the square of two reflects the natural derivation of the result.

The *viṣamakarman*, in turn, may be seen by considering that, in Fig. 1, the entire square has side $a$, and that this side is decomposed as $(a-b) + b$. The difference $a^2 - b^2$ is then equal to the square $(a-b)^2$ and two oblong quadrilaterals of area $ab$. Moving one of them, we reconstruct an oblong of area $(a+b)(a-b)$.

---

[49] A. Keller's translation (*Expounding the Mathematical Seed : Bhāskara I on the Mathematical Chapter of the Āryabhaṭīya,* two volumes, Birkhäuser, Berlin, 2006) appears to be the only one that literally translates *dvikṛti* (vol. 1, p. 103). She does not investigate its mathematical significance, because Bhāskara I does not, but correctly relates it to the *saṅkramaṇa* and its applications (vol. 2, p. 106 and vol. 1, p. xxxiv, n. 84).

[50] Bhāskara I similarly uses *guṇita* for *hata* in his commentary on II.7, and *kṛta* for four in his gloss on II.5, ex. 2.

[51] We write *g* for *guru* and *l* for *laghu* ("heavy" and "light" syllables respectively). These basic syllables may be grouped into *gaṇa*-s (groups). It is admissible to replace a group *llll* by, say, *gll* as here. For the rules that an āryā verse must satisfy, see the Appendix.

[52] Brahmagupta gives the same result in different words in the problem section at the end of chapter 18:
 18.100: *śeṣavadhād dvikṛtiguṇāt śeṣāntaravargasaṃyutānmūlam |*
   *śeṣāntaronayuktaṃ dalitaṃ śeṣe pṛthagabhīṣṭe ‖*
Even though Brahmagupta rephrased the result, he too expresses four as the square of two (*dvikṛti*), confirming that this formulation is significant.

It therefore appears that Āryabhaṭa knew that the *saṅkramaṇa* and the *viṣamakarman* were standard tools – he does not care to recall the former even though he uses it –, and sought to integrate them in a connected discourse that would give in one blow: the expansion of $(a+b)^2$, the relation $a^2 - b^2 = (a+b)(a-b)$, and $(a-b)^2 + 4ab = (a+b)^2$, from which the expansion of $(a-b)^2$ immediately follows as well. This line of argument is suggested by writing 4 as $2^2$, prompting the reader to scrutinize the text.

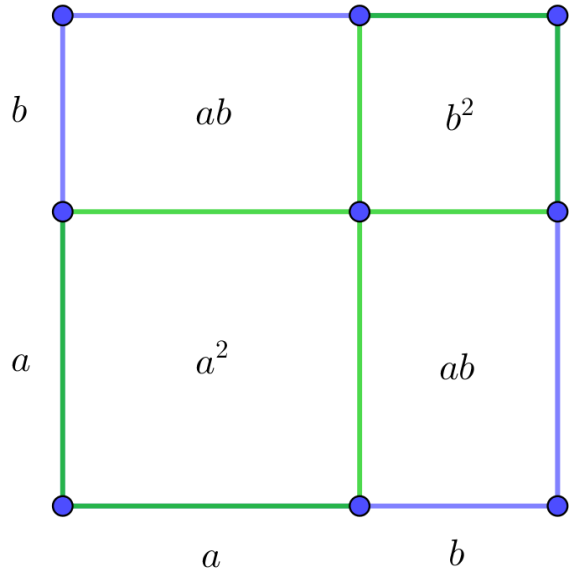

*Figure 1. Āryabhṭīya* II.23. The square of the sum $a+b$ is the complete square. It is decomposed into two squares and two (equal) oblongs. Removing the sum of the squares from the square of the sum as Āryabhaṭa prescribes, we are left with two oblongs of sides $a$ and $b$. This provides a geometric derivation of Prop. II.23.

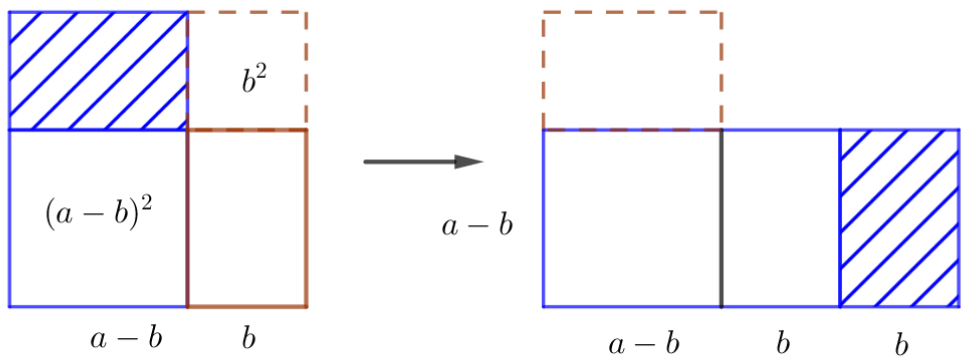

*Figure 2. Āryabhaṭīya* II.24. On the left, a square of area is removed from the entire square (of area $a^2$), leaving an area $a^2 - b^2$ represented by one square of area $(a-b)^2$ and two oblongs, both of area $ab$. The hatched area is then moved as indicated on the right, thus completing an oblong of dimensions $a-b$ and $a+b$. This provides a geometric derivation of the identity $a^2 - b^2 = (a+b)(a-b)$ on which Prop. II.24 rests.

### 5. Brahmagupta's *varga-prakṛti* and quadratic surds.

Starting from Prop. 18.65 in BSS, Brahmagupta introduces the theory of equation $b^2 - Na^2 = k$, called *varga-prakṛti*[53]. The numbers $N$ and $k$ are given, and $a$ and $b$ are sought. Brahmagupta's focus is on integral solutions, but rational solutions may come in as intermediate steps. This is one of the important topics of (Indian) Mathematics, and its modern theory still is not quite complete even though much has been accomplished[54]. We focus here on its relation to the irrationality of surds[55]. We show that Brahmagupta's treatment of *varga-prakṛti* implies a derivation of irrationality of $\sqrt{N}$ whenever the *varga-prakṛti* equation for one nonzero value of $k$ (for example, $k = 1$) has infinitely solutions. This will follow from a close reading of Prop. 18.65-66 and 18.73 in *BSS*. Even for $\sqrt{2}$, this seems to be an argument for irrationality that is not found anywhere else.

Recall that Brahmagupta had earlier described how to extend the usual operations with increasing generality: first to negative numbers (cf. *BSS* 18.31-36), after which he extended the *saṃkramaṇa* and *viṣamakarman* (18.37). He then introduced quadratic surds (*BSS* 18.38-41) to "unmanifest" letters (that is, he introduced a literal algebra in several variables) (*BSS* 18.42-43). This is followed by his theory of equations (18.44-64), and his treatment of *varga-prakṛti*. Now, while extending operations to expressions of the form $u + \sqrt{v}$, he restricted the word *karaṇī* mean only irrational square roots[56], since otherwise the term $\sqrt{v}$ could be merged with $u$. Brahmagupta also stressed in this connection that these surds are constructible[57]. However, he did not give any reason for the change, even though this move implies the recognition that some square roots are irrational, a statement made in no earlier Indian text. Such a departure from tradition must have been made for good reason, and one would expect this reason to be apparent in his text. It is not surprising that Brahmagupta does not state irrationality directly: it appears that Brahmagupta, like most Indian authors, *refrains from making negative statements* except in polemical discourse. This feature becomes less conspicuous in later texts, but even there, when we do find statements equivalent to impossibility results, the expression is softened. For instance, in Bhāskarācārya [58]'s criticisms against incorrect understandings [59] of Brahmagupta's area formula (Prop. 12.21, see *Līlāvatī*, 168-171) they are mostly phrased as rhetorical questions[60], and include a slur against the "ancients", which is a serious anomaly in an Indian context. It is apparent that criticisms must be softened if the opponent belongs to the same group as the author. Among scholars, criticism should remain civilized and therefore, remain implied rather than explicit as much as possible.

---

[53] This name should be preferred to "Pell's equation" or similar denominations.
[54] The modern theory is formulated in terms of properties of quadratic extensions of the field of rational numbers. One point of entry into this literature from a historical point of view is H. Hasse, "History of Class Field Theory", in *Algebraic Number Theory*, J.W.S. Cassels and A. Fröhlich (ed.), Academic Press, New York, 1967, pp. 266–279.
[55] For an overview of sources and results, see B.B. Datta and A.N. Singh, *History of Hindu Mathematics : A Source Book,* Part II, Asia Publishing House, Bombay, 1938, pp. 141-181. Recall that a surd is a root that cannot be expressed rationally. It is derived from Greek *alogos* 'inexpressible, irrational' through the Arabic (Oxford Eng. Dict.). The Latin *surdus* has several meanings, including "inaudible" and "deaf"; the former is relevant here. The term is usually applied to irrational *radicals* only.
[56] In the *Śulvasūtra*s, this word refers to the "maker", the cord that "makes" a square of which it is the side. There is no irrationality requirement there: the *karaṇī* of four square units is the cord of two units.
[57] On this process of successive generalizations, see our analysis of *BSS* 18.37, 18.38 and 18.42-43 in *Gaṇita Bhāratī* **41**, No. 1-2 (2019), 93-113, §5.2.
[58] Also known as Bhāskara II to distinguish him from the commentator of the same name.
[59] Bhāskara II does not realize that these were incorrect because at that time, the meaning of Brahmagupta's neologism *tricaturbhuja* had already been lost.

Let us turn to BSS 18.65-66 and follow Brahmagupta's argument.

18.65. *mūlaṃ dvidheṣṭavargād guṇakaguṇādiṣṭayutavihīnācca |*
*ādyavadho guṇakaguṇaḥ sahāntyaghātena kṛtamantyam ||*

The '[final] root' [is obtained from] the root of an arbitrary, squared,
Multiplied by the multiplier and increased or diminished by an arbitrary.
[A new] 'final [root]' is made up from the product of the [old] 'final' [roots] combined
With the product of the [old] 'initials' multiplied by the multiplier.

18.66. *vajravadhaikyam prathamam prakṣepaḥ kṣepavadhatulyaḥ |*
*prakṣepaśodhakahṛte mūle prakṣepake rūpe ||*

The union of the cross-products [of the roots] is the 'first'.
The [new] interpolator is equal to the product of the interpolators
One divides, by the [additive] interpolator or subtractive [interpolator],
Both roots when the [desired] interpolator is unity.

Here, Brahmagupta defines his problem and notation, and gives his composition rule in only two verses. Let us see how by scrutinizing his formulation. Brahmagupta introduces two sets of four variables. The text comprises three parts that will be analysed in turn: 18.65a, 18.65b-18.66a and 18.66b[61].

**18.65a.** In the first half-verse, we are directed to start from two 'roots' $b$ and $b'$ obtained by performing the same operation twice (*dvidhā*): take an arbitrary square $a^2$ (*iṣṭavarga*), multiply it by a multiplier $N$, and add or subtract from it another arbitrary quantity $k$; the root of the result is $b$. We will refer to $k$ as the interpolator[62] or the additive (resp. subtractive), as convenient. Similarly, $b'$ is obtained from another square $a'^2$, with a possibly different interpolator $\pm k'$, additive or subtractive. We therefore have

$$Na^2 \pm k = b^2 \text{ and } Na'^2 \pm k' = b'^2. \tag{VP}$$

---

[60] Such as: "Since the diagonals of the quadrilateral are indeterminate, how should the area be in this case determinate? The diagonals, found as assumed by the ancients, do not answer in another case. […] How can a person, neither specifying one of the perpendiculars, nor either of the diagonals, ask the rest? or how can he demand a determinate area, while they are indefinite?" (Colebrooke's translations.) The text reads: *caturbhujasyāniyatau hi karṇau kathaṃ tatao'sminniyatam phalam syāt | prasādhitau tacchravaṇau yadādyaiḥ sva-kalpitau tāvitaratra na staḥ | teṣuvevabāhuṣvaparau ca karṇāvanekadhā kṣetraphalam tatas ca || […] lambayoḥ karṇayorvaikam-anirdiśyāparaṃ katham | pṛcchatyaniyatatve'pi niyatañcāpi tatphalam ||.*
[61] We denote by a and b the halves of each verse.
[62] *Lit*.: "that which is thrown in", interpolated between $Na^2$ and $b^2$. It can be of any sign. Here, *kṣepa* may refer to the either of the two arbitrary quantities in the second quarter of 18.65 (cf. *iṣṭayutavihīnāt*), while *prakṣepa* denotes the new additive. When the two signs are to be stressed, Brahmagupta writes *prakṣepa-śodhaka.*

By saying *dvidhā* without further specification, the case $a' = a$, $b' = b$, $k = k'$ is not excluded.

Thus, for each problem, there are two given quantities: a 'multiplier' (*guṇaka*) $N$ and an 'arbitrary' (*iṣṭa*) $k$, as well as two 'roots' (*mūla*), the initial (*ādya*) and the final (*antya*), called here $a$ and $b$ respectively. Here, $N$ is a coefficient common to all the equations in this proposition, because it is not 'arbitrary'. These four quantities are related by the two equations (VP) described in the first hemistich.

**18.65b-18.66a.** The second hemistich of 18.65 and the first of 18.66 explain the generation of a new solution from two given ones. The new 'final' root is the "product of the initials times the multiplier, together with the product of the finals" or, in symbols, $b'' = Naa' + bb'$. 18.66a gives the new initial root, here called 'first' (*prathama*), possibly to clarify that it is neither of the 'initial' roots previously considered. It is the "union of the cross-products" of the roots, namely $a'' = ab' + a'b$. One says that the new solutions are obtained by composition[63] of $a, b$ with $a', b'$.

This result is found in most later works and is usually given a simple verification based on the expansion of a square. Recall that the expansion of a square is a consequence of the *saṅkramamaṇa* and *viṣamakarman*. The calculation is then immediate: Assume we have two sets of 'roots' satisfying the same equation: $b^2 - Na^2 = k$ ; $b'^2 - Na'^2 = k'$. We then have $Na''^2 = N(a^2b'^2 + a'^2b^2 + 2aa'bb')$, while $b''^2 = N^2a^2a'^2 + b^2b'^2 + 2Naa'bb'$, so that the combination $b''^2 - Na''^2$ equals $N^2a^2a'^2 + b^2b'^2 - N(a^2b'^2 + a'^2b^2) = (b^2 - Na^2)(b'^2 - Na'^2) = kk'$. Therefore, $a'', b''$ do provide a solution with interpolator $kk'$ as announced.

However, the origin of these formulae is not made clear by this calculation. This origin is probably to be sought in the theory of quadratic surds that Brahmagupta developed earlier in the same chapter. He dealt with quantities of the form $\sqrt{u} + v$. Quantities of the form $a\sqrt{N}$ were reduced to the standard form $\sqrt{Na^2}$. It follows from the rules for manipulation of surds that

$$\left(\sqrt{Na^2} + b\right)\left(\sqrt{Na^2} - b\right) = Na^2 - b^2 = -k. \tag{1}$$

Remember that Brahmagupta introduced manipulation rules for negative numbers as well as surds. Similarly,

$$\left(\sqrt{Na'^2} + b'\right)\left(\sqrt{Na'^2} - b'\right) = Na'^2 - b'^2 = -k'. \tag{2}$$

and

$$\left(\sqrt{Na''^2} + b''\right)\left(\sqrt{Na''^2} - b''\right) = Na''^2 - b''^2. \tag{3}$$

Now, if $N$ is positive[64], we have

---

[63] The usual Sanskrit term is *bhāvanā*.
[64] Brahmagupta does not consider square roots of negative numbers.

$$\left(\sqrt{Na^2}+b\right)\left(\sqrt{Na'^2}+b'\right) = b'\sqrt{Na^2}+b\sqrt{Na'^2}+Naa'+bb'$$
$$= \sqrt{Na^2b'^2}+\sqrt{Na'^2b^2}+Naa'+bb'. \qquad (4)$$

This suggests appealing to Proposition 18.38[65], that enables one to add surds: if $uu' = z^2$, then $\sqrt{u}+\sqrt{u'} = \sqrt{u+2z+u'}$. Since $(Na^2b'^2)(Na'^2b^2) = (Nab'a'b)^2$, we find from 18.38 that $\sqrt{Na^2b'^2}+\sqrt{Na'^2b^2}$ is the square root of $N(ab')^2 + 2N(ab')(a'b) + N(a'b)^2 = N(ab'+a'b)^2$. Substituting into (4), we have

$$\left(\sqrt{Na^2}+b\right)\left(\sqrt{Na'^2}+b'\right) = \sqrt{N(ab'+a'b)^2}+Naa'+bb' = \sqrt{Na''^2}+b'', \qquad (5)$$

with the values of $a''$ and $b''$ given in 18.65-66. Similarly,

$$\left(\sqrt{Na^2}-b\right)\left(\sqrt{Na'^2}-b'\right) = -\sqrt{N(ab'+a'b)^2}+Naa'+bb' = b''-\sqrt{Na''^2}. \qquad (6)$$

Multiplying (5) and (6) together, and using relations (1–3), we obtain:

$$\begin{aligned} Na''^2 - b''^2 &= \left(\sqrt{Na''^2}+b''\right)\left(\sqrt{Na''^2}-b''\right) \\ &= -\left(\sqrt{Na^2}+b\right)\left(\sqrt{Na'^2}+b'\right)\left(\sqrt{Na^2}-b\right)\left(\sqrt{Na'^2}-b'\right) \\ &= -\left(\sqrt{Na^2}+b\right)\left(\sqrt{Na^2}-b\right)\left(\sqrt{Na'^2}+b'\right)\left(\sqrt{Na'^2}-b'\right) \\ &= -(Na^2-b^2)(Na'^2-b'^2) = -kk'. \end{aligned}$$

Hence the relation

$$Na''^2 + kk' = b''^2. \qquad (7)$$

with the values of $a''$ and $b''$ given in 18.65-66. Thus, they are a natural outgrowth of Brahmagutpa's earlier results. However, the direct verification based on the expansion of the square of a sum has a distinct advantage: *it also applies if the multiplier N is negative.* It is likely that this traditional verification also goes back to Brahmagupta. It remains to deal with the second hemistich of 18.66.

**18.66b.** Here the singular *prakṣepaśodhaka* implies that the 'interpolators' are equal, and may be additive or subtractive. Let us therefore write their common value $\pm k$ with $k$ positive. For both choices of signs, the product of the (equal) interpolators is $k^2$, so that (7) now reads

$$Na''^2 + k^2 = b''^2.$$

If we divide the solutions by $k$, letting $x = a''/k$ and $y = b''/k$, we obtain, as Brahmagupta wrote,

---

[65] See our analysis of BSS 18.38 in *Gaṇita Bhāratī* **41** (1) (2019), 93-113.

$$Nx^2 + 1 = y^2.$$

This stresses the special position of the case when $k = 1$, and also stresses the usefulness of composing rational solutions in order to get integral ones. This observation is relevant in Prop. 18.67-68, that show how to solve with $k = 1$, we know a solution with[66] $k = \pm 4$. Also, Brahmagupta also tells us that by composing any two integral (resp. rational) solution of the equation $Na^2 + 1 = b^2$, we obtain another integral (resp. rational) solution *of the same equation*. Indeed, if $k = k' = 1$, we also have $kk' = 1$.

Assume now that $N$ is a positive integer, $k = 1$, and exclude the solutions[67] $a = 0, b = \pm 1$. Also assume that we know one other integer solution. Then there are *arbitrarily large solutions of the same equation*. More generally, if we have an integral solution of $Na'^2 + k' = b'^2$, with $a' > 0, b' > 0$, then we may find a solution of the same equation with $a' \geq a + 1$. Indeed, we may take $a$ and $b$ to be positive, since they are nonzero and only enter the equation through their squares. In that case, $a'' = ab' + a'b$ and $b'' = Naa' + bb'$ are also positive. If, as we assume, the solutions $a, b, a', b'$ are all nonzero integers, we necessarily have $N, a, b, a', b' \geq 1$ and therefore, $a'' \geq a' + 1$. If $k' = 1$, it suffices to have one solution since, as noted above, it is permissible to take $a = a', b = b'$. Iterating the process, we may generate arbitrarily large solutions. This leads to the two results. First, *if equation $Nx^2 + 1 = y^2$ has an integral solution, it has infinitely many, arbitrarily large solutions* by repeated compositions of this solution with itself. Second, by applying composition with arbitrary $k'$ and with $k = 1$, *if equations $Nx^2 + k' = y^2$ and $Nx^2 + 1 = y^2$ each have a positive integer solution, then each of them has infinitely, arbitrarily large many solutions*.

This completes the analysis of 18.65-66.

In the next few propositions, Brahmagupta describes shortcuts to get at solutions with $k = 1$, and then turns to the case when $N$ is a perfect square and therefore, $\sqrt{N}$ is rational. His result is a very simple application of the *saṃkramaṇa*. We now see that his formulation, combined with his earlier results, implies that it is not possible to get a solution with $k = 1$. It is this result that, by turning it around, yields a criterion for the irrationality of $\sqrt{N}$ in other cases.

> 18, 73. *varge guṇake kṣepaḥ kenaciduddhṛtayutonito dalitaḥ |*
>    *prathamo'ntyamūlamanyo guṇakārapadoddhṛtaḥ prathamaḥ ||*

> When the multiplier is a square, the interpolator
> Is divided by any [quantity] whatsoever, and added or subtracted, and the result halved.
> The first [of these results] is a 'final root'; the other,
> Divided by the square root of the multiplier, is a 'first [root]'.

---

[66] If we have a solution with $k = -1$, we obtain one with $k = 1$ by composing it with itself. Similarly, a solution with $k = \pm 2$, yields one with $k = 4$ by composing it with itself.

[67] Since $Na^2 + 1 > 0$, we cannot have $b = 0$.

Brahmagupta assumes that the multiplier is a perfect square: $N = n^2$. The text expresses that for any $m$ whatsoever,

$$a = \frac{1}{2n}\left(\frac{k}{m} - m\right), \quad b = \frac{1}{2}\left(\frac{k}{m} + m\right),$$

yields a solution of $n^2 a^2 + k = b^2$. This solution is easy to derive using the *saṃkramaṇa*. Indeed, the equation may be written $k = b^2 - (na)^2 = (b - na)(b + na)$, so that if we let $b - na = m$, we have $b + na = k/m$. The values of $b$ and $na$ are determined by their sum and difference, leading to the above solution. Thus, Brahmagupta's solution is not only valid, but is also the most general solution. Here, $n$ or $m$ could be rational or integral.

Given what Brahmagupta has expounded before, it is difficult to avoid asking whether there are integer solutions if $k = 1$. Brahmagupta's tools imply very simply that there is none if $N$ is the square of a rational number. Indeed, assume that $N = p^2/q^2$, with $p/q$ irreducible and $Na^2 + k = b^2$ with $k$ a positive integer. We may also assume that the integers $p, q, a, b$ are all positive. We then have

$$(pa)^2 + kq^2 = (bq)^2,$$

hence $(bq + pa)(bq - pa) = kq^2$ is independent of the solution $(a, b)$. This means that $bq + pa$ and $bq - pa$ are both bounded above: since an integer cannot be smaller than 1, $bq + pa$ and $bq - pa$ cannot exceed $kq^2$. But we have seen that if there is a solution with $k = 1$, there are necessarily arbitrarily large ones. Since $bq + pa$ is larger than both $a$ and $b$, $bq + pa$ also can be arbitrarily be large. This contradicts the above bound. Therefore, if $Na^2 + k = b^2$ has an integer solution, then $\sqrt{N}$ is irrational, it is a *karaṇī* in Brahmagupta's sense.

For instance, since $2a^2 + 1 = b^2$ has the solution $a = 2, b = 3$, $\sqrt{2}$ is irrational. Similarly, since $3a^2 + 1 = b^2$ has the solution $a = 1, b = 2$, $\sqrt{3}$ is irrational. However, $4a^2 + 1 = b^2$ has no solution in positive integers, however large. Indeed, we would need to have $b^2 > 4a^2$ hence $b > 2a$, and also $(b + 2a)(b - 2a) = 1$, But the product of two positive integers is equal to one if and only if both are equal to 1. We therefore obtain $b + 2a = b - 2a = 1$, hence $a = 0$, which contradicts the assumption.[68]

Brahmagupta avoids a negative statement by suggesting the reader to investigate whether composition theorems 18.65-66 apply here, that is, to investigate the coherence of his discourse.

### 6. Conclusions.

In this first part, we first recalled the difference between apodictic and dogmatic discourse and showed that apodictic discourse – a conclusive and motivated discourse meant to be analyzed rather than interpreted –, is quite common in ancient and contemporary Mathematics, particularly in research papers. We then obtained the following results.

---

[68] As an exercise, the reader may want to work out the details for $5a^2 + 1 = b^2$ (with the solution $a = 4, b = 9$), $6a^2 + 1 = b^2$ (with the solution $a = 2, b = 5$), and $9a^2 + 1 = b^2$ (no positive solution).

- In India, because of the avoidance of writing, the development of a discursive view of science encouraged the replacement of sense objects and actions by word-representations. This led in turn to a second level of abstraction in which word-representations were themselves manipulated by an intellectual process, rather than subjected to a formal process that could be carried out in writing. This accounts for the abstract and general nature of propositions and the very precise terminology found in Indian mathematical texts.
- The term *vyavahāra* for Brahmagupta, and probably for most later mathematical authors, refers to an apodictic discourse when it is applied to a part of a mathematical treatise. This is consistent the meaning of this word in other contexts.
- Two further examples of apodictic discourses have been examined and have led to the following conclusions.
    - Āryabhaṭa's treatment of the expansion of a square in Prop. II.23 of his *Āryabhaṭīya* was meant to suggest a geometric derivation.
    - Āryabhaṭa suggested a particular chain of argument by his seemingly quaint replacement of 4 by the square of 2 in his Prop. II.24.
    - Brahmagupta related the irrationality of surds to the existence of arbitrarily large solutions of the *varga-prakṛti* problem through the contrast between Prop. 18.65-66 and 18.73 of his *Brāhmasphuṭasiddhānta*.

### Appendix : the *āryā* metre

Since *BSS* and the *Āryāṣṭaśata* – as its name indicates –, are in *āryā*[69] metre, we reproduce here its description[70] from Piṅgala's *Chandassūtra*. Recall that Sanskrit syllables are distinguished into *laghu* ('light', denoted by *l*) and *guru* ('heavy', denoted by *g*)[71]. The *ja-gaṇa* (denoted by *j*) is *lgl*, and the *na-gaṇa* is *lll*. Therefore, a *nagaṇa-la* consists of four *laghus* (*llll*). Here is the text and its translation.

IV.14 *svarārdhaṃ cāryā''rdhaṃ*
IV.15 *atrā'yuṅ na j*
IV.16 *ṣaṣṭho j*
IV.17 *nlau vā*
IV.18 *nlau cet padaṃ dvitīyādi*
IV.19 *saptamaḥ prathamādi*
IV.20 *antye pañcamaḥ*
IV.21 *ṣaṣṭhaś ca l*

IV.14: And half an *āryā* is seven[72] and a half (groups).
IV.15: there, no *jagaṇa* in odd (groups),

---

[69] It hardly needs to be recalled that the word *āryā* has no ethnic undertones in an Indian context.
[70] We followed A. Weber ("Über die Metrik der Inder. Zwei Abhandlungen", *Indische Studien* VIII, Harrwitz und Goßmann, Berlin, 1863, se p. 291). See also L. Renou et J. Filliozat, *L'Inde Classique : Manuel des Études Indiennes*, Maisonneuve, Paris, 1985, Appendice 2, or the appendix on verses in Apte's Dictionary that has an extensive list of metres.
[71] The guru and laghu are also often represented by an avagraha and a daṇḍa respectively.
[72] Seven is indicated by *svara* because there are seven notes in the octave.

IV.16: (and) the sixth (group) is a *jagaṇa*,
IV.17: or a *nagaṇa-la*.
IV.18: If a *nagaṇa-la*, a word starts at the second (*laghu*);
IV.19: if (it is) the seventh (group), (a word) starts at the first (*laghu* of these four);
IV.20: (if it is) the fifth (group) in the final (half, a word starts) from the first light.
IV.21: And the sixth (group of this final half) is (an isolated) *laghu*.

In addition, a caesura is optional after the third group, in each hemistich. It is very often found in *BSS*. This caesura divides each half-verse into four quarters, which is why we translated as often as possible *āryā* verses into four lines in English. If both hemistichs have this caesura, the verse is called *pathyā*, otherwise it is *vipulā* (IV.22–23). The standard form has *lgl* or *llll* as sixth *gaṇa* of the first hemistich, and *l* as sixth *gaṇa* of the second hemistich. If both halves follow the pattern of the first half, it is called a *gīti* (this is the metre of Āryabhaṭa's *Daśagītikā*); if they both follow the pattern of the second half, it is an *upagīti*; if the two patterns are interchanged, we have an *udgīti*; and if there are eight full *gaṇa*s in each half, it is an *āryāgīti*.[73]

This meter used to be very popular – this is the metre of the *Sāṃkhyakārikā*, also translated by Colebrooke[74]. Metres based on *gaṇa*s are not found in the Vedas, and their origin is disputed. The *āryā* is found in Prākṛt and Marāṭhī poetry[75]. The *āryā* metre is quite restrictive and difficult because there are restrictions on the structure of all four quarter-verses – not just on the total number of syllables, or on the structure of two of the four *pada*s only –, unlike other common metres. Its knowledge is often useful in reading mathematical texts. For instance, it is the metre in *Ābh* II.1 that implies that its author's name is *Āryabhaṭa* rather than *Āryabhaṭṭa*. Indeed, the second hemistich begins with *āryabhaṭastviha* (*gll gll*) whereas *āryabhaṭṭastviha* would give *glg gll* instead, with a faulty first *gaṇa*.[76]

**Contact details:** Satyanad Kichenassamy, Professor of Mathematics,
Université de Reims Champagne-Ardenne,
Laboratoire de Mathématiques de Reims (CNRS, UMR9008), B.P. 1039, F-51687 Reims Cedex 2.
*E-mail* : satyanad.kichenassamy@univ-reims.fr
*Web* : https://www.normalesup.org/~kichenassamy/

---

[73] *Chandaḥsūtra*, IV.28-31: *ādyaradhasamā gītiḥ | antyenopasgītiḥ | utkrameṇodgītiḥ | ardhe vasugaṇa āryāgītiḥ* (see Weber, *op. cit.,* pp. 302-303).
[74] *The Sánkhya Káriká, or Memorial Verses on the Sánkhya Philosophy, by Íswara Krishna; translated from the Sanskrit by Henry Thomas Colebrooke, Esq., and also the Bháshya of Commentary by Gaurapáda* [*sic*]*; translated, and illustrated by an original comment by Horace Hyman Wilson, M.A. F.R.S.* (Oxford, London, 1837).
[75] This innovation is perhaps also to be compared to the notion of *cīr* in Tamil prosody, that has no Sanskrit analogue, since any *gaṇa* with four *mātrā*s is equivalent to a pair of *cīr*.
[76] See also *Ābh* I.1 and IV.50.